\documentclass[letterpaper, 10 pt, conference]{ieeeconf}  

\IEEEoverridecommandlockouts                              
\usepackage{epsfig} 
\usepackage{epstopdf}

\usepackage{fancyhdr}
\usepackage{indentfirst}
\usepackage{amsmath}
\usepackage{amssymb}
\usepackage{graphicx}
\usepackage{color}
\usepackage{latexsym}
\usepackage{amsfonts}
\usepackage{amssymb,amsmath,amsfonts}
\usepackage{amsmath,mathrsfs,bm,url,times}
\usepackage{cite}
\usepackage{caption}
\usepackage{subcaption}

\newtheorem{theorem}{Theorem}
\newtheorem{lemma}{Lemma}

\newtheorem{definition}{Definition}

\newtheorem{remark}{Remark}

\DeclareMathOperator{\argmax}{argmax}
\DeclareMathOperator{\diag}{diag}
\title{\LARGE \bf
Distributed Dynamic Event-Triggered Control for Multi-Agent Systems
}

\author{Xinlei Yi, Kun Liu, Dimos V. Dimarogonas and Karl H. Johansson
\thanks{This work was supported by the
Knut and Alice Wallenberg Foundation, the  Swedish Foundation for Strategic Research, the Swedish Research Council, and the National Natural Sciences Foundation of China under Grant 61503026.}
\thanks{X. L. Yi, D. V. Dimarogonas and K. H. Johansson are with the ACCESS Linnaeus Centre, Electrical Engineering, KTH Royal Institute of Technology, 100 44, Stockholm, Sweden; K. Liu is with School of Automation, Beijing Institute of Technology, 100081 Beijing, China
.
        {\tt\small \{xinleiy, dimos, kallej\}@kth.se, kunliubit@bit.edu.cn}.}%
}

\begin{document}

\maketitle
\thispagestyle{empty}
\pagestyle{empty}

\begin{abstract}\label{s:Abstract}We propose two distributed dynamic triggering laws to solve the consensus problem for multi-agent systems with event-triggered control. Compared with existing triggering laws, the proposed triggering laws involve internal dynamic variables which play an essential role to guarantee that the triggering time sequence does not exhibit Zeno behavior. Some existing triggering laws are special cases of our dynamic triggering laws. Under the condition that the underlying graph is undirected and connected, it is proven that  the proposed dynamic triggering laws together with the event-triggered control make the state of each agent converges exponentially to the average of the agents' initial states. Numerical simulations illustrate the effectiveness of the theoretical results and show that the dynamic triggering laws lead to reduction of actuation updates and inter-agent communications.
\end{abstract}
\section{INTRODUCTION}\label{sec:intro}

Multi-agent (average) consensus problem, where a group of agents seeks to agree upon certain quantity of interest (e.g., the average of their initial states), has been widely investigated because it has many applications such as mobile robots, autonomous underwater vehicles, unmanned air vehicles, etc. There are many results obtained in
this field, such as \cite{olfati2004consensus,ren2007information,Liu2011consensus} and the references therein. In these papers, agents have continuous-time dynamics and actuation. However, in practice, it is  in most cases at discrete points in time that agents communicate with their neighbors and take action. There are also many papers that study the agents with discrete-time dynamics or continuous-time dynamics with discontinuous information
transmission, for example see \cite{xiao2008asynchronous,liu2010necessary,you2011network}. In these papers, time-driven sampling is used to determine when agents should establish communication with its neighbors. Time-driven sampling is often implemented by periodic sampling. A significant drawback of periodic sampling is that it requires all agents to exchange their information synchronously, which is not so easy to be realized in real systems, especially when the number of agents is large.

In addition to time-driven sampling, event-driven sampling has been proposed \cite{aastrom1999comparison,tabuada2007event}. In event-driven sampling actuation updates and inter-agent communications occur only when some specific events are triggered, for instance, a measure of the state error exceeds a specified threshold. Event-driven sampling is normally implemented by event-triggered or self-triggered control. The event-triggered control is often piecewise constant between
the triggering times. The triggering times are determined by the triggering laws. Many researchers studied event-triggered control for multi-agent systems recently  \cite{dimarogonas2012distributed,garcia2013decentralised,seyboth2013event,nowzari2016distributed,meng2013event,fan2013distributed,yi2016distributed,yi2016formation,yi2017pull,fan2015self}.
A key challenge in event-triggered control for multi-agent systems is how to design triggering laws to determine the corresponding triggering times, and to exclude Zeno behavior. For continuous-time multi-agent systems, Zeno behavior means that there are infinite number of triggers in a
finite time interval \cite{johansson1999regularization}. Another important question is how to realize the event-triggered controller in a distributed way.

In \cite{girard2015dynamic}, by introducing an internal dynamic variable, a new class of event-triggering mechanisms is presented. The idea of using internal dynamic variables in event-triggered and self-triggered control can also be found in \cite{dolk2015dynamic,de2013robust}. In this paper, we modify the dynamic event triggering mechanism in \cite{girard2015dynamic} and extend it to multi-agent systems in a distributed manner.

We have the following main contributions: we propose two dynamic triggering laws  which are distributed in the sense that they do not require any a priori knowledge of global network parameters; we prove that the proposed dynamic triggering laws yield consensus exponentially fast;  and we show that they are free from Zeno behavior. We show also that the triggering laws in \cite{dimarogonas2012distributed,garcia2013decentralised,seyboth2013event} are special cases of the control laws considered in this paper.

The rest of this paper is organized as follows. \mbox{Section \ref{secpreliminaries}} introduces the preliminaries and the problem formulation. The main results are stated in Section \ref{secmain}. Simulations are given in Section \ref{secsimulation}. Finally, the paper is concluded in Section \ref{secconclusion}.

\noindent {\bf Notations}: $\|\cdot\|$ represents the Euclidean norm for
vectors or the induced 2-norm for matrices. ${\bf 1}_n$ denotes the column
vector with each component being 1 and dimension $n$. $I_n$ is the $n$ dimension identity matrix. $\rho_2(\cdot)$ indicates the minimum
positive eigenvalue for matrices having positive eigenvalues. Given two
symmetric matrices $M,N$, $M\ge N$ means $M-N$ is a positive semi-definite matrix. $|S|$ is the cardinality of set $S$.

\section{PRELIMINARIES}\label{secpreliminaries}
In this section, we present some definitions from algebraic graph theory \cite{mesbahi2010graph} and the formulation of the problem.

\subsection{Algebraic Graph Theory}

Let $\mathcal G=(\mathcal V,\mathcal E, A)$ denote a (weighted) undirected graph with the set of
agents (vertices or nodes) $\mathcal V =\{1,\dots,n\}$, the set of links (edges) $\mathcal E
\subseteq \mathcal V \times \mathcal V$, and the (weighted) adjacency matrix
$A =A^{\top}=(a_{ij})$ with nonnegative elements $a_{ij}$. A link
 $(i,j)\in \mathcal E$ if $a_{ij}=a_{ji}>0$, i.e., if agent $i$ and $j$ can communicate with each other. It is assumed that $a_{ii}=0$ for all $i\in \mathcal V$. Let $\mathcal{N}_i=\{j\in \mathcal V\mid a_{ij}>0\}$ and $\deg_i=\sum\limits_{j=1}^{n}a_{ij}$ denotes the neighbors and (weighted) degree of agent $i$, respectively. The degree matrix of graph $\mathcal G$ is $D=\diag(\deg_1, \cdots, \deg_n)$. The Laplacian matrix is $L=(L_{ij})=D-A$. A  path of length $k$ between agent $i$ and agent $j$ is a subgraph with distinct agents $i_0=i,\dots,i_k=j\in\mathcal V$ and edges $(i_j,i_{j+1})\in\mathcal E,~j=0,\dots,k-1$. a subgraph with distinct agents $i_0=i,\dots,i_k=j\in\mathcal V$ and edges $(i_j,i_{j+1}),~j=0,\dots,k-1$.

\begin{definition}
An undirected graph is  connected if there exists at least one path between any two agents. And an undirected graph is completed if any two distinct agents are connected by an edge.
\end{definition}

Obviously, there is a one-to-one correspondence between a graph and its adjacency matrix or its Laplacian matrix. If we let $K_n=I_n-\frac{1}{n}\bf{1}_n\bf{1}_n^{\top}$, then we can treat $K_n$ as the Laplacian matrix of a completed graph with $n$ agents and edge weight $\frac{1}{n}$.

For a connected graph we have the following well known results.
\begin{lemma}
If a graph $\mathcal G$ is connected, then its Laplacian matrix $L$ is positive semi-definite, i.e., $z^{\top}Lz\ge0$ for any $z\in\mathbb{R}^n$. And, $z^{\top}Lz=0$ if and only if $z=a\bf{1}_n$ for some $a\in\mathbb R$. Moreover, we have
\begin{align}\label{LKn}
0\le\rho_2(L)K_n\le L.
\end{align}
\end{lemma}
{\bf Proof}: For the proof of (\ref{LKn}), please see Lemma 2.1 in \cite{yi2016formation}.

\subsection{Problem Formulation}
We consider a set of $n$ agents that are modelled as a single integrator
\begin{align}
\dot{x}_i(t)=u_i(t)
,~i\in\mathcal V, t\ge0,\label{system}
\end{align}
where $x_i(t)\in\mathbb{R}$ is the state and $u_i(t)\in\mathbb{R}$ is the control input.
\begin{remark}
For the ease of presentation, we study the case where all the agents have scalar states. However, the analysis in this paper is also valid for the cases where the agents have vector-valued states.
\end{remark}

In the literature, the following distributed consensus protocol is often considered, e.g., \cite{olfati2004consensus,ren2007information},
\begin{align}\label{inputc}
u_i(t)=-\sum_{j=1}^{n}L_{ij}x_j(t).
\end{align}

To implement the consensus control protocol (\ref{inputc}), continuous-time state information from neighbours is needed. However, it is often impractical to require continuous communication in physical applications.

Inspired by the idea of event-triggered control for multi-agent systems \cite{dimarogonas2012distributed}, we use instead of (\ref{inputc}) the following event-triggered control
\begin{align}\label{input}
u_i(t)=-\sum_{j=1}^{n}L_{ij}x_j(t^j_{k_j(t)}),
\end{align}
where $k_{j}(t)=\argmax_{k}\{t^{j}_{k}\le t\}$ with the increasing  $\{t_{k}^{j}\}_{k=1}^{\infty}$, $j\in\mathcal V$ to be determined later. We assume $t^j_1=0, j\in\mathcal V$.
Note that the control protocol (\ref{input}) only updates at the triggering times and is constant between consecutive triggering times.

For simplicity, let $x(t)=[x_{1}(t),\dots,x_{n}(t)]^{\top}$, $\hat{x}_{i}(t)=x_{i}(t_{k_{i}(t)}^{i})$,
$\hat{x}(t)=[\hat{x}_{1}(t), \dots,\hat{x}_{n}(t)]^{\top}$,  $ e_{i}(t)=\hat{x}_{i}(t)-x_{i}(t)$, and $e(t)=[e_{1}(t),\cdots,e_{n}(t)]^{\top}=\hat{x}(t)-x(t)$. Then we can rewrite system (\ref{system}) with even-triggered control (\ref{input}) in the following stack vector form:
\begin{align}
\dot{x}(t)=-L\hat{x}(t)=-L(x(t)+e(t)).
\end{align}

\section{DYNAMIC EVENT-TRIGGERED CONTROL}\label{secmain}
In this section, we propose the dynamic triggering laws to determine the triggering time sequence  and we prove that they lead to consensus for the multi-agent system (\ref{system}) with event-triggered control (\ref{input}).
\subsection{Continuous Approach}
We first give the following lemma.
\begin{lemma}\label{lemma2}
Consider the multi-agent system (\ref{system}) with event-triggered control (\ref{input}).
The average of all agents' states $\bar{x}(t)=\frac{1}{n}\sum_{i=1}^{n}x_i(t)$ is a constant, i.e., $\bar{x}(t)=\bar{x}(0),\forall t\ge0$.
\end{lemma}
{\bf Proof}: It follows from (\ref{system}) and (\ref{input}) that the time derivative of the average value is given by
\begin{align*}
\dot{\bar{x}}(t)=&\frac{1}{n}\sum_{i=1}^{n}\dot{x}_i(t)=-\frac{1}{n}\sum_{i=1}^{n}\sum_{j=1}^{n}L_{ij}x_j(t^j_{k_j(t)})\\
=&-\frac{1}{n}\sum_{j=1}^{n}x_j(t^j_{k_j(t)})\sum_{i=1}^{n}L_{ij}=0.
\end{align*}
Thus $\bar{x}(t)$ is a constant.

Consider a Lyapunov candidate:
\begin{align}\label{V}
V(t)=&\frac{1}{2}x^{\top}(t)K_nx(t)=\frac{1}{2}x^{\top}(t)[I_n-\frac{1}{n}{\bf1}_n{\bf1}_n^{\top}] x(t)\nonumber\\
=&\frac{1}{2}\sum_{i=1}^{n}x_i^2(t)-\frac{n}{2}\bar{x}^2(0)=\frac{1}{2}\sum_{i=1}^{n}[x_i(t)-\bar{x}(0)]^2.
\end{align}
Then the derivative of $V(t)$ along  the trajectories of system (\ref{system}) with the event-triggered control (\ref{input}) satisfies
\begin{align}
\dot{V}(t)
=&\sum_{i=1}^{n}[x_i(t)-\bar{x}(0)]\dot{x}_i(t)\nonumber\\
=&\sum_{i=1}^{n}x_i(t)\dot{x}_i(t)-\bar{x}(0)\sum_{i=1}^{n}\dot{x}_i(t)
=\sum_{i=1}^{n}x_i(t)\dot{x}_i(t)\nonumber\\
=&\sum_{i=1}^{n}x_i(t)\sum_{j=1}^{n}(-L_{ij}x_j(t^j_{k_j(t)}))\nonumber
\end{align}
\begin{align}
=&-\sum_{i=1}^{n}x_i(t)\sum_{j=1}^{n}L_{ij}(x_j(t)+e_j(t))\nonumber\\
\overset{*}{=}&-\sum_{i=1}^{n}q_i(t)-\sum_{i=1}^{n}\sum_{j=1}^{n}x_i(t)L_{ij}e_j(t)\nonumber\\
=&-\sum_{i=1}^{n}q_i(t)-\sum_{i=1}^{n}\sum_{j=1}^{n}e_i(t)L_{ij}x_j(t)\nonumber\\
=&-\sum_{i=1}^{n}q_i(t)-\sum_{i=1}^{n}\sum_{j=1,j\neq i}^{n}e_i(t)L_{ij}(x_j(t)-x_i(t))\nonumber\\
\le&-\sum_{i=1}^{n}q_i(t)-\sum_{i=1}^{n}\sum_{j=1,j\neq i}^{n}L_{ij}e_i^2(t)\nonumber\\
&-\sum_{i=1}^{n}\sum_{j=1,j\neq i}^{n}L_{ij}\frac{1}{4}(x_j(t)-x_i(t))^2\nonumber\\
=&-\sum_{i=1}^{n}q_i(t)+\sum_{i=1}^{n}L_{ii}e_i^2(t)\nonumber\\
&-\sum_{i=1}^{n}\sum_{j=1}^{n}\frac{1}{4}L_{ij}(x_j(t)-x_i(t))^2\nonumber\\
\overset{*}{=}&-\sum_{i=1}^{n}\frac{1}{2}q_i(t)+\sum_{i=1}^{n}L_{ii}e_i^2(t),\label{dV}
\end{align}
where
\begin{align}\label{qi}
q_{i}(t)=-\frac{1}{2}\sum_{j=1}^{n}L_{ij}(x_{j}(t)-x_{i}(t))^2\ge0,
\end{align}
and the equalities denoted by $\overset{*}{=}$ hold since
\begin{align*}
\sum_{i=1}^{n}q_{i}(t)
=&-\sum_{i=1}^{n}\frac{1}{2}\sum_{j=1}^{n}L_{ij}(x_{j}(t)-x_{i}(t))^2\\
=&\sum_{i=1}^{n}\sum_{j=1}^{n}x_{i}(t)L_{ij}x_{j}(t)=x^{\top}(t)Lx(t),
\end{align*}
and the inequality holds since $ab\le a^2+\frac{1}{4}b^2$.

Similar to \cite{dimarogonas2012distributed} and \cite{yi2016distributed}, if we use the following law to determine the triggering time sequence:
\begin{align}
t^i_1=&0\nonumber\\
t^i_{k+1}=&\max_{r\ge t^i_k}\Big\{r:e^2_i(t)\le \frac{\sigma_i}{2L_{ii}}q_i(t),
\forall t\in[t^i_k,r]\Big\},\label{statictriggersingle}
\end{align}
with $\sigma_i\in(0,1)$,
then, from (\ref{dV}) and (\ref{statictriggersingle}), we have
\begin{align}
\dot{V}(t)
&\le-\sum_{i=1}^{n}\frac{1}{2}q_i(t)+\sum_{i=1}^{n}L_{ii}e^2_i(t)\nonumber\\
&\le-\frac{1}{2}(1-\sigma_{\max})\sum_{i=1}^{n}q_i(t)\nonumber\\
&=-\frac{1}{2}(1-\sigma_{\max})x^{\top}(t)Lx(t)\nonumber\\
&\le-\frac{1}{2}(1-\sigma_{\max})\rho_2(L)x^{\top}(t)K_nx(t)\nonumber\\
&=-(1-\sigma_{\max})\rho_2(L)V(t),\label{dVstatic}
\end{align}
where $\sigma_{\max}=\max\{\sigma_1,\dots,\sigma_n\}<1$ and the last inequality holds due to (\ref{LKn}).
Then
\begin{align}\label{Vupps}
V(t)\le V(0)e^{-(1-\sigma_{\max})\rho_2(L)t}.
\end{align}
This implies that system (\ref{system}) reaches consensus exponentially.

\begin{remark}\label{staticremark}
We call (\ref{statictriggersingle})  a static triggering law since it does not involve any  extra dynamic variables except $x_i(t),\hat{x}_i(t)$ and $x_j(t),j\in\mathcal{N}_i$. The static triggering law (\ref{statictriggersingle}) is  distributed since each agent's control action only depends on its neighbours' state information, without any prior knowledge of global parameters, such as the eigenvalue of the Laplacian matrix.
\end{remark}

\begin{remark}\label{dimos}
If we consider the same graph that considered in \cite{dimarogonas2012distributed}, i.e., $a_{ij}=1$ if $(i,j)\in \mathcal E$, then $L_{ii}=|N_i|$. From the facts $a(1-a|N_i|)\le\frac{1}{4|N_i|}$ and $(\sum_{j=1}^{n}(x_j(t)-x_i(t)))^2\le2|N_i|\sum_{j=1}^{n}(x_j(t)-x_i(t))^2$, we have
$\frac{\sigma_ia(1-a|N_i|)}{|N_i|}(\sum_{j=1}^{n}(x_j(t)-x_i(t)))^2\le\frac{\sigma_i}{2|N_i|}q_i(t)$.
In other words, the distributed triggering law (10) in \cite{dimarogonas2012distributed} is a special case of the static triggering law (\ref{statictriggersingle}).
\end{remark}

The main purpose of using the event-triggered control is to reduce the overall need of  actuation updates and communication between agents, so it is essential to exclude Zeno behavior. However, we do not know whether Zeno behavior can be excluded or not in the above event-triggered control law. In order to explicitly exclude Zeno behavior, in the following we propose a dynamic event-triggered control law.

Inspired by \cite{girard2015dynamic}, we propose the following internal dynamic variable $\eta_i$ to agent $i$:
\begin{align}\label{etai}
\dot{\eta}_i(t)=-\beta_i\eta_i(t)+\xi_i(\frac{\sigma_i}{2}q_i(t)-L_{ii}e^2_i(t)), i\in\mathcal V,
\end{align}
with $\eta_i(0)>0$, $\beta_i>0$, $\xi_i\in[0,1]$, and $\sigma_i\in[0,1)$. These dynamic variables are correlated in the event-triggered law, as defined in our first main result.
\begin{theorem}\label{dynamictheorem}
Consider the multi-agent system (\ref{system}) with the event-triggered control protocol (\ref{input}).
Suppose that the underlying graph $\mathcal G$ is undirected and connected.
Given $\theta_i>\frac{1-\xi_i}{\beta_i}$ and the first triggering time $t^i_1=0$, agent $i$ determines the triggering time sequence $\{t^i_k\}_{k=2}^{\infty}$ by
\begin{align}\label{dynamictriggersingle}
t^i_{k+1}=\max_{r\ge t^i_k}\Big\{r:\theta_i\Big(L_{ii}e^2_i(t)-\frac{\sigma_i}{2}q_i(t)\Big)
\le \eta_i(t),&\nonumber\\
\forall t\in[t^i_k,r]&\Big\},
\end{align}
with $q_i(t)$ defined in (\ref{qi}) and $\eta_i(t)$ defined in (\ref{etai}).
Then the consensus is achieved exponentially and there is no Zeno behavior.
\end{theorem}
{\bf Proof}: (i) From equation (\ref{etai}) and condition (\ref{dynamictriggersingle}), we have
$$\dot{\eta}_i(t)\ge -\beta_i\eta_i(t)-\frac{\xi_i}{\theta_i}\eta_i(t).$$
Thus
\begin{align}\label{etailow}
\eta_i(t)\ge\eta_i(0)e^{-(\beta_i+\frac{\xi_i}{\theta_i})t}>0.
\end{align}

Consider a Lyapunov candidate:
\begin{align}\label{W}
W(t)=V(t)+\sum_{i=1}^{n}\eta_i(t).
\end{align}
Then the derivative of $W(t)$ along the trajectories of systems (\ref{system}) and  (\ref{etai}) with the event-triggered control (\ref{input}) satisfies
\begin{align}
&\dot{W}(t)=\dot{V}(t)+\sum_{i=1}^{n}\dot{\eta}_i(t)\nonumber\\
\le&-\sum_{i=1}^{n}\frac{1}{2}q_i(t)+\sum_{i=1}^{n}L_{ii}e^2_i(t)-\sum_{i=1}^{n}\beta_i\eta_i(t)\nonumber\\
&+\sum_{i=1}^{n}\xi_i(\frac{\sigma_i}{2}q_i(t)-L_{ii}e^2_i(t))\nonumber\\
=&-\sum_{i=1}^{n}\frac{1}{2}(1-\sigma_i)q_i(t)-\sum_{i=1}^{n}\beta_i\eta_i(t)\nonumber\\
&+\sum_{i=1}^{n}(\xi_i-1)(\frac{\sigma_i}{2}q_i(t)-L_{ii}e^2_i(t))\nonumber\\
\le&-\sum_{i=1}^{n}\frac{1}{2}(1-\sigma_i)q_i(t)-\sum_{i=1}^{n}\beta_i\eta_i(t)
+\sum_{i=1}^{n}\frac{1-\xi_i}{\theta_i}\eta_i(t)\nonumber\\
=&-\sum_{i=1}^{n}\frac{1}{2}(1-\sigma_i)q_i(t)-\sum_{i=1}^{n}\Big(\beta_i-\frac{1-\xi_i}{\theta_i}\Big)\eta_i(t)\nonumber\\
\le&-(1-\sigma_{\max})\sum_{i=1}^{n}\frac{1}{2}q_i(t)-k_d\sum_{i=1}^{n}\eta_i(t)\nonumber\\
\le&-(1-\sigma_{\max})\rho_2(L)V(t)-k_d\sum_{i=1}^{n}\eta_i(t)\nonumber\\
\le&-k_WW(t),
\end{align}
where $k_d=\min_{i}\{\beta_i-\frac{1-\xi_i}{\theta_i}\}>0$ and $k_W=\min\{(1-\sigma_{\max})\rho_2(L),k_d\}>0$. Then
\begin{align}
V(t)<W(t)\le W(0)e^{-k_Wt}.\label{Vupp}
\end{align}
This implies that system (\ref{system}) reaches consensus exponentially.

(ii) Next, we prove that there is no Zeno behavior by contradiction. Suppose there exists Zeno behavior. Then there exists an agent $i$, such that $\lim_{k\rightarrow+\infty}t^i_k=T_0$ where $T_0$ is a positive constant.

From (\ref{Vupp}), we know that there exists a positive constant $M_0>0$ such that $|x_i(t)|\le M_0$ for all $t\ge0$ and $i=1,\dots,n$.
Then, we have $$|u_i(t)|\le2M_0L_{ii}.$$

Let $\varepsilon_0=\frac{\eta_i(0)}{4\sqrt{\theta_iL_{ii}^3}M_0}e^{-\frac{1}{2}(\beta_i+\frac{1}{\theta_i})T_0}>0$. Then from the property of limit,  there exists a positive integer $N(\varepsilon_0)$ such that
\begin{align}\label{zeno}
t^i_k\in[T_0-\varepsilon_0,T_0],~\forall k\ge N(\varepsilon_0).
\end{align}
Noting $q_i(t)\ge0$ and (\ref{etailow}), we can conclude that one sufficient condition to guarantee the inequality in condition (\ref{dynamictriggersingle}) is
\begin{align}\label{suffi1}
|\hat{x}_i(t)-x_i(t)|\le\frac{\eta_i(0)}{\sqrt{\theta_iL_{ii}}}e^{-\frac{1}{2}(\beta_i+\frac{\xi_i}{\theta_i})t}.
\end{align}
Again noting $|\dot{x}_i(t)|=|u_i(t)|\le2M_0L_{ii}$ and $|\hat{x}_i(t^i_{k})-x_i(t^i_{k})|=0$ for any triggering time $t^i_k$, we can conclude that one sufficient condition to the above inequality is
\begin{align}\label{suffi2}
(t-t^i_k)2M_0L_{ii}\le\frac{\eta_i(0)}{\sqrt{\theta_iL_{ii}}}e^{-\frac{1}{2}(\beta_i+\frac{\xi_i}{\theta_i})t}.
\end{align}
Then
\begin{align}
&t^i_{N(\varepsilon_0)+1}-t^i_{N(\varepsilon_0)}
\ge\frac{\eta_i(0)}{2\sqrt{\theta_iL_{ii}^3}M_0}e^{-\frac{1}{2}(\beta_i+\frac{\xi_i}{\theta_i})t^i_{N(\varepsilon_0)+1}}\nonumber\\
&\ge\frac{\eta_i(0)}{2\sqrt{\theta_iL_{ii}^3}M_0}e^{-\frac{1}{2}(\beta_i+\frac{\xi_i}{\theta_i})T_0}=2\varepsilon_0,
\end{align}
which contradicts to (\ref{zeno}). Therefore, Zeno behavior is excluded.

\begin{remark}\label{dynamicremark}
We call (\ref{dynamictriggersingle}) a dynamic triggering law since it involves the extra dynamic variables $\eta_i(t)$. Similar to the static triggering law (\ref{statictriggersingle}), it is also distributed.
The static triggering law (\ref{statictriggersingle}) can be seen as a limit case of the dynamic triggering  law (\ref{dynamictriggersingle}) when $\theta_i$ grows large. Thus, from the analysis in Remark \ref{dimos}, we can conclude that the distributed triggering law (9) in \cite{dimarogonas2012distributed} is a special case of the dynamic triggering law (\ref{dynamictriggersingle}).
\end{remark}

\begin{remark}\label{seyboth}
If we choose $\xi_i=0$ in (\ref{etai}) and $\sigma_i=0$ in (\ref{dynamictriggersingle}), then $\eta_i(t)=\eta_i(0)e^{-\beta_it}$ and now the inequality in (\ref{dynamictriggersingle}) is $|e_i(t)|\le\frac{\sqrt{\eta_i(0)}}{\sqrt{\theta_iL_{ii}}}e^{-\frac{\beta_i}{2}t}$. This is the triggering  function (7) in \cite{seyboth2013event} with $c_0=0,c_1=\frac{\sqrt{\eta_i(0)}}{\sqrt{\theta_iL_{ii}}}, \alpha=\frac{\beta_i}{2}$. However, we do not need the constraint $\alpha<\rho_2(L)$ which is necessary in \cite{seyboth2013event}.
\end{remark}


If we choose $\beta_i$ large enough, then $k_W=(1-\sigma_{\max})\rho_2(L)$. Hence, in this case, from (\ref{Vupps}) and (\ref{Vupp}), we know that the trajectories of the multi-agent system (\ref{system}) with the event-triggered control (\ref{input}) under static event-triggered control law (\ref{statictriggersingle}) and dynamic event-triggered control law (\ref{dynamictriggersingle}) have the same guaranteed decay rate given by (\ref{Vupps}).

\subsection{Discontinuous Approach}
In the above static and dynamic triggering control laws, in order to check the inequalities (\ref{statictriggersingle}) and (\ref{dynamictriggersingle}), each agent still needs to continuously monitor its neighbors's states, which means continuous communication is still needed. In the following, we will modify the above results to avoid this.

We upper-bound the derivative of $V(t)$ along  the trajectories of system (\ref{system}) with the event-triggered control (\ref{input}) by a different way.
Similar to the derivation process to get (\ref{dV}), we have
\begin{align}
\dot{V}(t)
=&\sum_{i=1}^{n}x_i(t)\sum_{j=1}^{n}-L_{ij}\hat{x}_j(t)\nonumber\\
=&-\sum_{i=1}^{n}(\hat{x}_i(t)-e_i(t))\sum_{j=1}^{n}L_{ij}\hat{x}_j(t)\nonumber\\
\overset{**}{=}&-\sum_{i=1}^{n}\hat{q}_i(t)+\sum_{i=1}^{n}\sum_{j=1}^{n}e_i(t)L_{ij}\hat{x}_j(t)\nonumber
\end{align}
\begin{align}
=&-\sum_{i=1}^{n}\hat{q}_i(t)+\sum_{i=1}^{n}\sum_{j=1,j\neq i}^{n}e_i(t)L_{ij}(\hat{x}_j(t)-\hat{x}_i(t))\nonumber\\
\le&-\sum_{i=1}^{n}\hat{q}_i(t)-\sum_{i=1}^{n}\sum_{j=1,j\neq i}^{n}L_{ij}e_i^2(t)\nonumber\\
&-\sum_{i=1}^{n}\sum_{j=1,j\neq i}^{n}L_{ij}\frac{1}{4}(\hat{x}_j(t)-\hat{x}_i(t))^2\nonumber\\
=&-\sum_{i=1}^{n}\hat{q}_i(t)+\sum_{i=1}^{n}L_{ii}e_i^2(t)\nonumber\\
&-\sum_{i=1}^{n}\sum_{j=1}^{n}\frac{1}{4}L_{ij}(\hat{x}_j(t)-\hat{x}_i(t))^2\nonumber\\
\overset{**}{=}&-\sum_{i=1}^{n}\frac{1}{2}\hat{q}_i(t)+\sum_{i=1}^{n}L_{ii}e_i^2(t),\label{dVdis}
\end{align}
where
\begin{align}\label{qih}
\hat{q}_{i}(t)=-\frac{1}{2}\sum_{j=1}^{n}L_{ij}(\hat{x}_{j}(t)-\hat{x}_{i}(t))^2\ge0,
\end{align}
and the equalities denoted by $\overset{**}{=}$ hold since
\begin{align*}
\sum_{i=1}^{n}\hat{q}_{i}(t)
=&-\sum_{i=1}^{n}\frac{1}{2}\sum_{j=1}^{n}L_{ij}(\hat{x}_{j}(t)-\hat{x}_{i}(t))^2\\
=&\sum_{i=1}^{n}\sum_{j=1}^{n}\hat{x}_{i}(t)L_{ij}\hat{x}_{j}(t)=\hat{x}^{\top}(t)L\hat{x}(t),
\end{align*}
and the inequality holds since $ab\le a^2+\frac{1}{4}b^2$.

Similar to \cite{garcia2013decentralised} and \cite{yi2016distributed}, if we use the following law to determine the triggering time sequence:
\begin{align}
t^i_1=&0\nonumber\\
t^i_{k+1}=&\max_{r\ge t^i_k}\Big\{r:e^2_i(t)\le \frac{\sigma_i}{2L_{ii}}\hat{q}_i(t),
\forall t\in[t^i_k,r]\Big\},\label{statictriggersingledis}
\end{align}
with $\sigma_i\in(0,1)$,
then, from (\ref{dVdis}) and (\ref{statictriggersingledis}), we have
\begin{align}
\dot{V}(t)
&\le-\sum_{i=1}^{n}\frac{1}{2}\hat{q}_i(t)+\sum_{i=1}^{n}L_{ii}e^2_i(t)\nonumber\\
&\le-\frac{1}{2}(1-\sigma_{\max})\sum_{i=1}^{n}\hat{q}_i(t)\nonumber\\
&=-\frac{1}{2}(1-\sigma_{\max})\hat{x}^{\top}(t)L\hat{x}(t).
\end{align}
Noting
\begin{align}
x^{\top}(t)Lx(t)
&=(\hat{x}(t)+e(t))^{\top}L(\hat{x}(t)+e(t))\nonumber\\
&\le2\hat{x}^{\top}(t)L\hat{x}(t)+2e^{\top}(t)Le(t)\nonumber\\
&\le2\hat{x}^{\top}(t)L\hat{x}(t)+2\|L\|\|e(t)\|^2\nonumber\\
&\le2\hat{x}^{\top}(t)L\hat{x}(t)
+\frac{\|L\|\sigma_{\max}}{\min_{i}L_{ii}}\sum_{i=1}^{n}\hat{q}_i(t)\nonumber\\
&=\Big(2+\frac{\|L\|\sigma_{\max}}{\min_{i}L_{ii}}\Big)
\hat{x}^{\top}(t)L\hat{x}(t),\label{xlx}
\end{align}
where the first inequality holds since $L$ is positive semi-definite and $a^{\top}Lb\le2a^{\top}La+2b^{\top}Lb,\forall a,b\in\mathbb{R}^n$, the second inequality holds since $a^{\top}La\le\|L\|\|a\|^2,\forall a\in\mathbb{R}^n$, and the last inequality holds due to (\ref{statictriggersingledis}),
we then have
\begin{align*}
\dot{V}(t)&\le-\frac{(1-\sigma_{\max})\min_{i}L_{ii}}{4\min_{i}L_{ii}+2\|L\|\sigma_{\max}}
x^{\top}(t)Lx(t)\\
&=-\frac{(1-\sigma_{\max})\min_{i}L_{ii}}{2\min_{i}L_{ii}+\|L\|\sigma_{\max}}\rho_2(L)x^{\top}(t)K_nx(t)\\
&=-\frac{(1-\sigma_{\max})\min_{i}L_{ii}}{2\min_{i}L_{ii}+\|L\|\sigma_{\max}}\rho_2(L)V(t).
\end{align*}
Then
\begin{align}\label{Vuppsdis}
V(t)\le V(0)e^{-\frac{(1-\sigma_{\max})\min_{i}L_{ii}}{2\min_{i}L_{ii}+\|L\|\sigma_{\max}}\rho_2(L)t}.
\end{align}
This implies that system (\ref{system}) reaches consensus exponentially.
\begin{remark}\label{garcia}
Similar to the analysis in Remark \ref{staticremark}, (\ref{statictriggersingledis}) is a static triggering law and it is also distributed. Moreover, similar to the analysis in Remark \ref{dimos}, we can conclude that the distributed  triggering law (6) in \cite{garcia2013decentralised} is a special case of the static triggering law (\ref{statictriggersingledis}).
\end{remark}

Just as the comment given in \cite{nowzari2016distributed} that the distributed  triggering law (6) in \cite{garcia2013decentralised} ``does not discard the possibility of an infinite number of events happening in a finite time period'', we also
do not know whether Zeno behavior can be excluded or not in the static triggering law (\ref{statictriggersingledis}). In the following, in order to explicitly exclude Zeno behavior, we will change the static triggering law (\ref{statictriggersingledis}) to the dynamic one.

Similar to (\ref{etai}), we propose an internal dynamic variable $\chi_i$ to agent $i$:
\begin{align}\label{chii}
\dot{\chi}_i(t)=-\beta_i\chi_i(t)+\xi_i(\frac{\sigma_i}{2}\hat{q}_i(t)-L_{ii}e^2_i(t)), i\in\mathcal V
\end{align}
with $\chi_i(0)>0$, $\beta_i>0$, $\xi_i\in[0,1]$, and $\sigma_i\in[0,1)$. Our second main result is given in the following theorem.
\begin{theorem}\label{dynamictheoremdis}
Consider the multi-agent system (\ref{system}) with the event-triggered control protocol (\ref{input}).
Suppose that the underlying graph $\mathcal G$ is undirected and connected.
Given $\theta_i>\frac{1-\xi_i}{\beta_i}$ and the first triggering time $t^i_1=0$, agent $i$ determines the triggering time sequence $\{t^i_k\}_{k=2}^{\infty}$ by
\begin{align}\label{dynamictriggersingledis}
t^i_{k+1}=\max_{r\ge t^i_k}\Big\{r:\theta_i\Big(L_{ii}e^2_i(t)-\frac{\sigma_i}{2}\hat{q}_i(t)\Big)
\le \chi_i(t),&\nonumber\\
\forall t\in[t^i_k,r]&\Big\},
\end{align}
with $\hat{q}_i(t)$ defined in (\ref{qih}) and $\chi_i(t)$ defined in (\ref{chii}).
Then the consensus is achieved exponentially and there is no Zeno behavior.
\end{theorem}
{\bf Proof}: (i) Similar to (\ref{etailow}), we have
\begin{align}
\chi_i(t)\ge\chi_i(0)e^{-(\beta_i+\frac{\xi_i}{\theta_i})t}>0.
\end{align}

Consider a Lyapunov candidate:
\begin{align}\label{F}
F(t)=V(t)+\sum_{i=1}^{n}\chi_i(t).
\end{align}
Then the derivative of $F(t)$ along the trajectories of systems (\ref{system}) and (\ref{chii}) with the event-triggered control (\ref{input}) satisfies
\begin{align}
&\dot{F}(t)=\dot{V}(t)+\sum_{i=1}^{n}\dot{\chi}_i(t)\nonumber\\
\le&-\sum_{i=1}^{n}\frac{1}{2}\hat{q}_i(t)+\sum_{i=1}^{n}L_{ii}e^2_i(t)-\sum_{i=1}^{n}\beta_i\chi_i(t)\nonumber\\
&+\sum_{i=1}^{n}\xi_i(\frac{\sigma_i}{2}\hat{q}_i(t)-L_{ii}e^2_i(t))\nonumber\\
=&-\sum_{i=1}^{n}\frac{1}{2}(1-\sigma_i)\hat{q}_i(t)-\sum_{i=1}^{n}\beta_i\chi_i(t)\nonumber\\
&+\sum_{i=1}^{n}(\xi_i-1)(\frac{\sigma_i}{2}\hat{q}_i(t)-L_{ii}e^2_i(t))\nonumber\\
\le&-\sum_{i=1}^{n}\frac{1}{2}(1-\sigma_i)\hat{q}_i(t)-\sum_{i=1}^{n}\beta_i\chi_i(t)
+\sum_{i=1}^{n}\frac{1-\xi_i}{\theta_i}\chi_i(t)\nonumber\\
=&-\sum_{i=1}^{n}\frac{1}{2}(1-\sigma_i)\hat{q}_i(t)-\sum_{i=1}^{n}\Big(\beta_i-\frac{1-\xi_i}{\theta_i}\Big)\chi_i(t)\nonumber\\
\le&-(1-\sigma_{\max})\sum_{i=1}^{n}\frac{1}{2}\hat{q}_i(t)-k_d\sum_{i=1}^{n}\chi_i(t)\nonumber\\
=&-\frac{1}{2}(1-\sigma_{\max})\hat{x}^{\top}(t)L\hat{x}(t)-k_d\sum_{i=1}^{n}\chi_i(t).
\end{align}
Similar to the derivation process to get (\ref{xlx}), we have
\begin{align}
&x^{\top}(t)Lx(t)\le2\hat{x}^{\top}(t)L\hat{x}(t)+2\|L\|\|e(t)\|^2\nonumber\\
\le&2\hat{x}^{\top}(t)L\hat{x}(t)
+\frac{\|L\|\sigma_{\max}}{\min_{i}L_{ii}}\sum_{i=1}^{n}\hat{q}_i(t)\nonumber\\
&+\frac{2\|L\|}{\min_{i}\{\theta_iL_{ii}\}}\sum_{i=1}^{n}\chi_i(t)\nonumber\\
=&\Big(2+\frac{\|L\|\sigma_{\max}}{\min_{i}L_{ii}}\Big)
\hat{x}^{\top}(t)L\hat{x}(t)+\frac{2\|L\|}{\min_{i}\{\theta_iL_{ii}\}}\sum_{i=1}^{n}\chi_i(t)\nonumber\\
\le & k_x\hat{x}^{\top}(t)L\hat{x}(t)+\frac{2\|L\|}{\min_{i}\{\theta_iL_{ii}\}}\sum_{i=1}^{n}\chi_i(t),\label{xlxdis}
\end{align}
where
\begin{align*}
k_x=\max\bigg\{2+\frac{\|L\|\sigma_{\max}}{\min_{i}L_{ii}},\frac{2(1-\sigma_{\max})\|L\|}{k_d\min_{i}\{\theta_iL_{ii}\}}\bigg\}.
\end{align*}
Then
\begin{align*}
&-\frac{1}{2}(1-\sigma_{\max})\hat{x}^{\top}(t)L\hat{x}(t)\\
\le&-\frac{1}{2k_x}(1-\sigma_{\max})x^{\top}(t)Lx(t)+\frac{k_d}{2}\sum_{i=1}^{n}\chi_i(t).
\end{align*}
Thus
\begin{align*}
\dot{F}(t)\le&-\frac{1}{2k_x}(1-\sigma_{\max})x^{\top}(t)Lx(t)-\frac{k_d}{2}\sum_{i=1}^{n}\chi_i(t)\\
\le&-\frac{\rho_{2}(L)}{2k_x}(1-\sigma_{\max})x^{\top}(t)K_nx(t)-\frac{k_d}{2}\sum_{i=1}^{n}\chi_i(t)\\
=&-\frac{\rho_{2}(L)}{k_x}(1-\sigma_{\max})V(t)-\frac{k_d}{2}\sum_{i=1}^{n}\chi_i(t)\\
\le& k_F F(t),
\end{align*}
where $k_F=\min\{\frac{\rho_{2}(L)}{k_x}(1-\sigma_{\max}),\frac{k_d}{2}\}$.
Then
\begin{align}
V(t)<F(t)\le F(0)e^{-k_Ft}.\label{Vuppdis}
\end{align}
This implies that system (\ref{system}) reaches consensus exponentially.

(ii) The way to exclude Zeno behavior is the same as the proof in Theorem \ref{dynamictheorem}.

\begin{remark}
Obviously, the triggering law (\ref{dynamictriggersingledis}) is dynamic and it is also distributed.
One can easily check that every agent does not need to continuously access its neighbors' states when  implementing the static and dynamic triggering laws (\ref{statictriggersingledis}) and (\ref{dynamictriggersingledis}).
\end{remark}

\begin{remark}
The static triggering law (\ref{statictriggersingledis}) can be seen as a limit case of the dynamic triggering law (\ref{dynamictriggersingledis}) when $\theta_i$ grows large. Thus, from the analysis in Remark \ref{garcia}, we can conclude that the distributed triggering law (6) in \cite{garcia2013decentralised} is a special case of the dynamic triggering tlaw (\ref{dynamictriggersingledis}).
\end{remark}

If we choose $\beta_i$ large enough, then $k_F=\frac{(1-\sigma_{\max})\min_{i}L_{ii}}{2\min_{i}L_{ii}+\|L\|\sigma_{\max}}\rho_2(L)$. Hence, in this case, from (\ref{Vuppsdis}) and (\ref{Vuppdis}), we know that the trajectories of the multi-agent system (\ref{system}) with event-triggered control (\ref{input}) under static event-triggered control law (\ref{statictriggersingledis}) and dynamic event-triggered control law (\ref{dynamictriggersingledis}) have the same guaranteed decay rate given by (\ref{Vuppsdis}).

\section{SIMULATIONS}\label{secsimulation}

\begin{figure}
\begin{subfigure}{.5\textwidth}
  \centering
  \includegraphics[width=.9\linewidth]{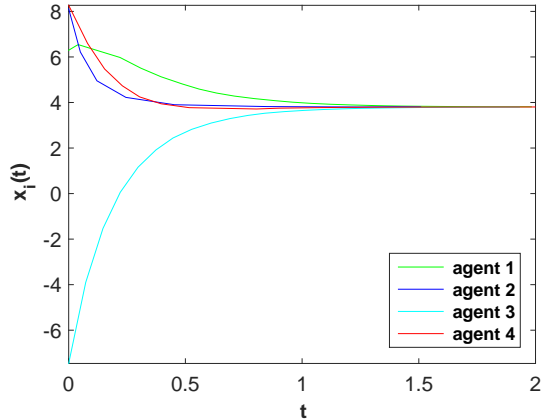}
  \caption{}
  \label{fig:2a}
\end{subfigure}%
\\
\begin{subfigure}{.5\textwidth}
  \centering
  \includegraphics[width=.9\linewidth]{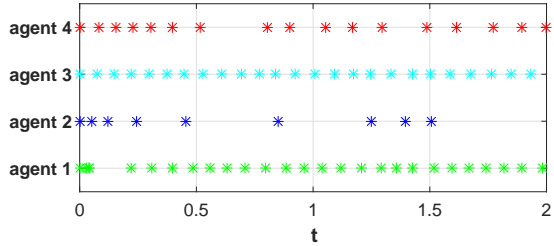}
  \caption{}
  \label{fig:2b}
\end{subfigure}
\caption{(a) The state evolution under the static triggering law (\ref{statictriggersingle}). (b) The triggering times for each agent under static triggering law (\ref{statictriggersingle}).}
\label{fig:2}
\end{figure}

\begin{figure}
\begin{subfigure}{.5\textwidth}
  \centering
  \includegraphics[width=.9\linewidth]{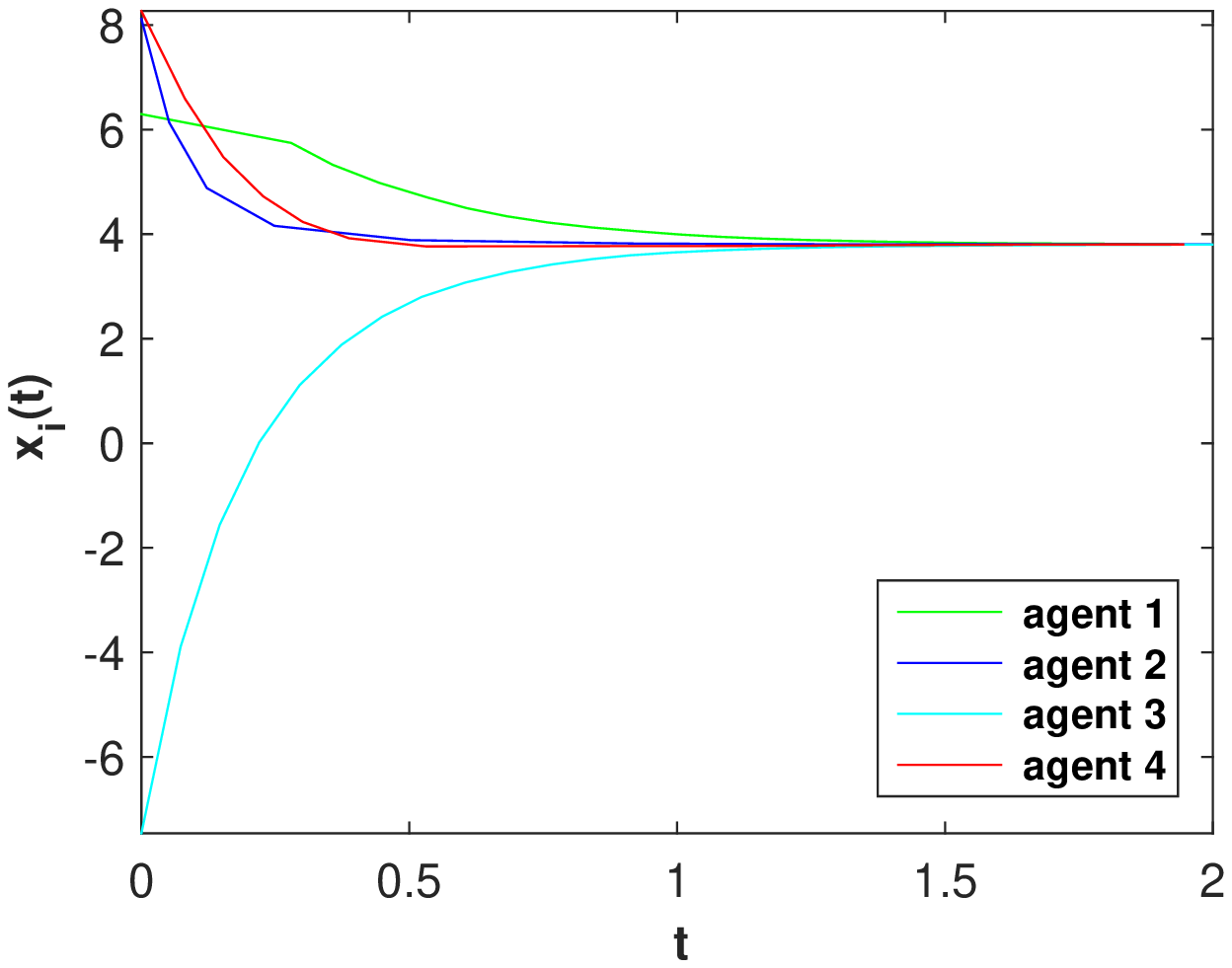}
  \caption{}
  \label{fig:3a}
\end{subfigure}%
\\
\begin{subfigure}{.5\textwidth}
  \centering
  \includegraphics[width=.9\linewidth]{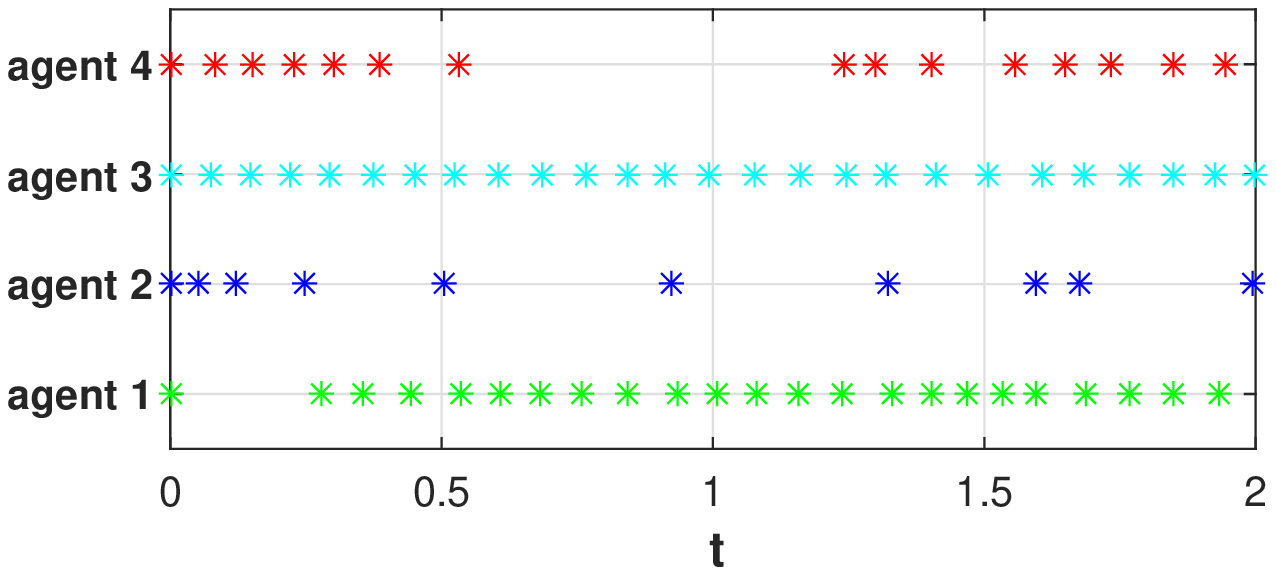}
  \caption{}
  \label{fig:3b}
\end{subfigure}
\caption{(a) The state evolution under the dynamic triggering law (\ref{dynamictriggersingle}). (b) The triggering times for each agent under dynamic triggering law (\ref{dynamictriggersingle}).}
\label{fig:3}
\end{figure}

\begin{figure}
\begin{subfigure}{.5\textwidth}
  \centering
  \includegraphics[width=.9\linewidth]{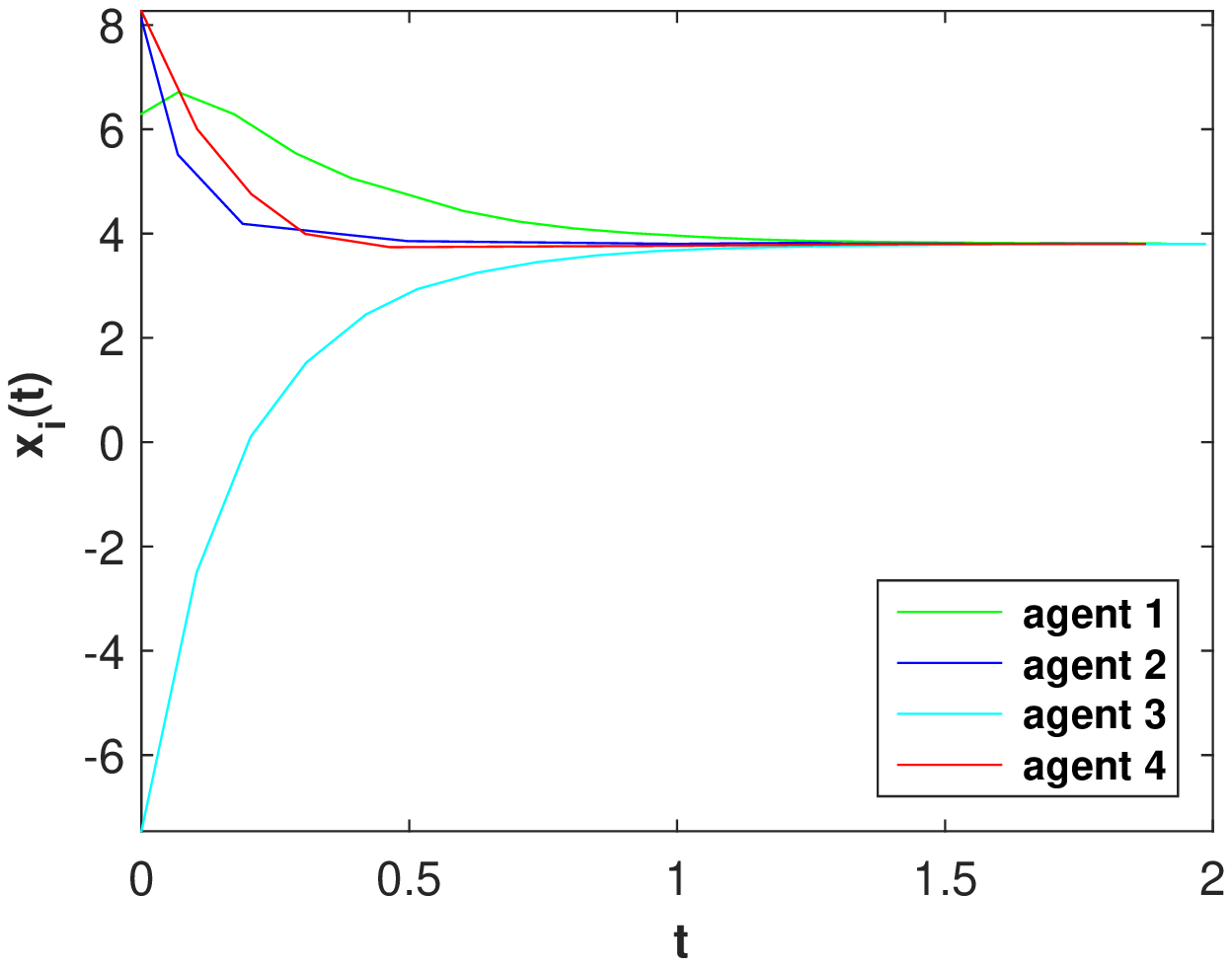}
  \caption{}
  \label{fig:4a}
\end{subfigure}%
\\
\begin{subfigure}{.5\textwidth}
  \centering
  \includegraphics[width=.9\linewidth]{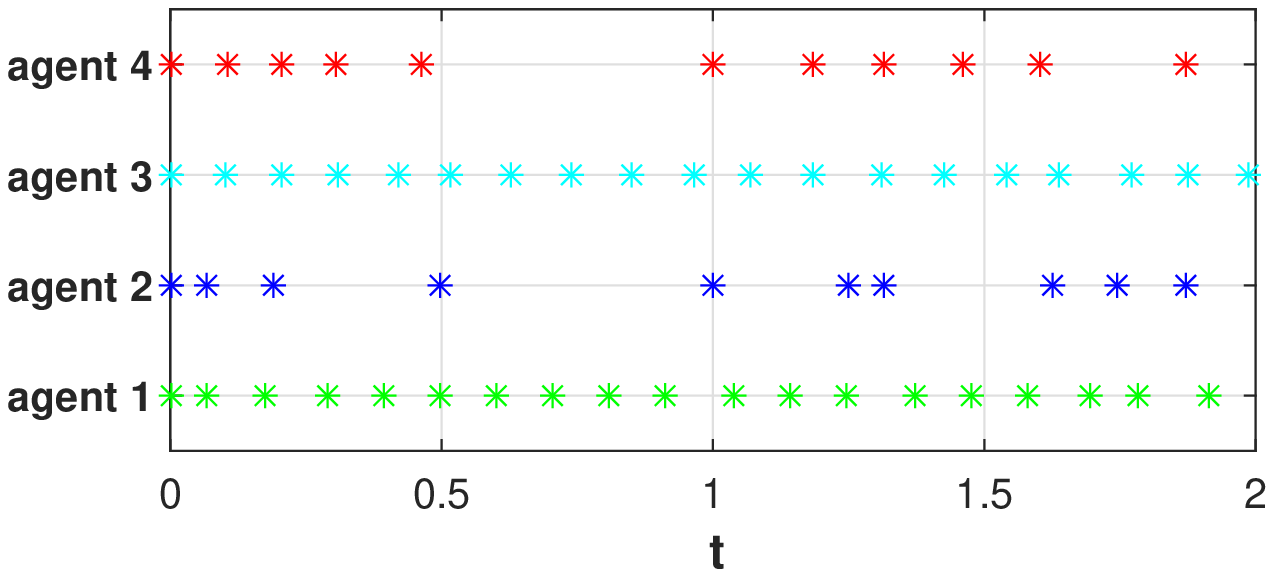}
  \caption{}
  \label{fig:4b}
\end{subfigure}
\caption{(a) The state evolution under the static triggering law (\ref{statictriggersingledis}). (b) The triggering times for each agent under static triggering law (\ref{statictriggersingledis}).}
\label{fig:4}
\end{figure}

\begin{figure}
\begin{subfigure}{.5\textwidth}
  \centering
  \includegraphics[width=.9\linewidth]{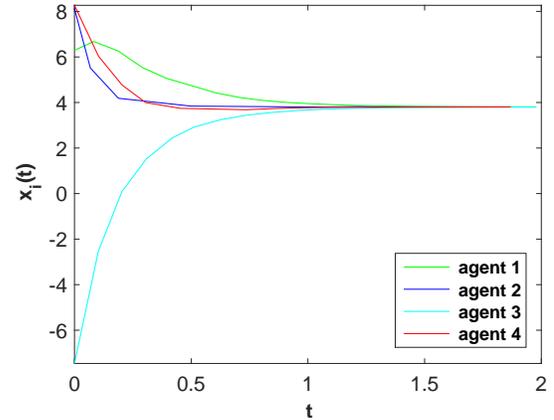}
  \caption{}
  \label{fig:5a}
\end{subfigure}%
\\
\begin{subfigure}{.5\textwidth}
  \centering
  \includegraphics[width=.9\linewidth]{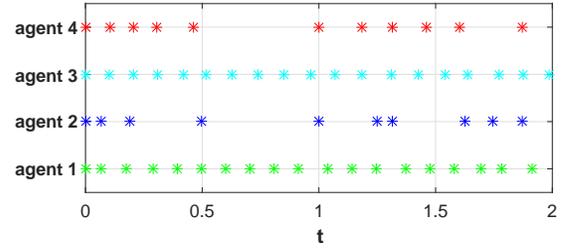}
  \caption{}
  \label{fig:5b}
\end{subfigure}
\caption{(a) The state evolution under the dynamic triggering law (\ref{dynamictriggersingledis}). (b) The triggering times for each agent under dynamic triggering law (\ref{dynamictriggersingledis}).}
\label{fig:5}
\end{figure}

In this section, a numerical example is given to demonstrate the presented results.
Consider a connected network of four agents with the Laplacian matrix
\begin{eqnarray*}
L=\left[\begin{array}{rrrr}3.4&-3.4&0&0\\
-3.4&9.8&-2.1&-4.3\\
0&-2.1&3.2&-1.1\\
0&-4.3&-1.1&5.4
\end{array}\right].
\end{eqnarray*}
The initial value of each agent is randomly selected within the interval $[-10,10]$. First, $x(0)=[6.2945, 8.1158, -7.4603, 8.2675]^{\top}$, the average initial state is $\bar{x}(0)=3.8044$. Fig. \ref{fig:2} (a) shows the state evolution under the static triggering law (\ref{statictriggersingle}) with $\sigma_i=0.5$. Fig. \ref{fig:2} (b) shows the corresponding triggering times for each agent. Fig. \ref{fig:3} (a) shows the state evolution under the dynamic triggering law (\ref{dynamictriggersingle}) with $\sigma_i=0.5$, $\eta_i(0)=10$, $\beta_i=1$, $\xi_i=1$ and $\theta_i=1$. Fig. \ref{fig:3} (b) shows the corresponding triggering times for each agent. Fig. \ref{fig:4} (a) shows the state evolution under the static triggering law (\ref{statictriggersingledis}) with $\sigma_i=0.5$. Fig. \ref{fig:4} (b) shows the corresponding triggering times for each agent. Fig. \ref{fig:5} (a) shows the state evolution under the dynamic triggering law (\ref{dynamictriggersingledis}) with $\sigma_i=0.5$, $\chi_i(0)=10$, $\beta_i=1$, $\xi_i=1$ and $\theta_i=1$. Fig. \ref{fig:5} (b) shows the corresponding triggering times for each agent.
It can
be seen that consensus is achieved when performing the four triggering laws proposed in this paper. Moreover, just as Theorem \ref{dynamictheorem} and Theorem \ref{dynamictheoremdis} point out, from the simulations we can also see that there is no Zeno behavior under  the dynamic triggering law (\ref{dynamictriggersingle}) and the dynamic triggering law (\ref{dynamictriggersingledis}). Although there is also no Zeno behavior under the static triggering law (\ref{statictriggersingle}) and the static triggering law (\ref{statictriggersingledis}) in the simulations, we still do not know how to prove this in theory.

\section{CONCLUSION}\label{secconclusion}
In this paper, we presented two dynamic triggering laws for multi-agent systems with event-triggered control.  We showed that,  some existing  triggering laws are special cases of the proposed dynamic triggering laws and  if the communication graph is undirected and connected, consensus is achieved exponentially. In addition, Zeno behavior was excluded by proving that the
triggering time sequence of each agent is divergent. Without any modifications, the results in this paper can be extended to the cases that the underlying graphs are directed, strongly connected and weight-balanced.  Future research directions
include considering general linear multi-agent systems and dynamic self-triggered control.

\bibliographystyle{IEEEtran}
\bibliography{refs}

\end{document}